\documentclass[12pt]{amsart}

\usepackage{enumerate}
\usepackage{amsmath, amsthm,amsfonts,amssymb,latexsym}
\usepackage[all]{xy}
\newtheorem{theorem}{Theorem}

\newtheorem{lemma}[theorem]{Lemma}

\newtheorem{corollary}[theorem]{Corollary}
\newtheorem{definition}[theorem]{Definition}

\newtheorem{remark}[theorem]{Remark}

\DeclareMathOperator{\Ran}{Ran}

\DeclareMathOperator{\Sym}{Sym}

\DeclareMathOperator{\co}{co}
\DeclareMathOperator{\rank}{rank}
 
\DeclareMathOperator{\Tip}{Tip}
\DeclareMathOperator{\Rep}{Rep}
\DeclareMathOperator{\Norm}{Norm}
\DeclareMathOperator{\Span}{Span}
\DeclareMathOperator{\NonTip}{NonTip}

\newcommand{\A}{\mathcal{A}}
\newcommand{\B}{\mathcal{B}}
\newcommand{\C}{\mathcal{C}}
\newcommand{\F}{\mathcal{F}}
\newcommand{\G}{\mathcal{G}}

\newcommand{\J}{\mathcal{J}}
\newcommand{\M}{\mathcal{M}}
\newcommand{\N}{\mathcal{N}}

\newcommand{\Real}{\mathbb{R}}
\newcommand{\Complex}{\mathbb{C}}

\newcommand{\abs}[1]{\left\vert#1\right\vert}
\newcommand{\set}[1]{\left\{#1\right\}}
\newcommand{\seq}[1]{\left<#1\right>}
\newcommand{\norm}[1]{\left\Vert#1\right\Vert}

\DeclareMathOperator{\diag}{diag}

\newcommand{\algebra}[1]{\Complex\langle #1 \rangle}

\begin{document}

\title[Flat functionals]{  Positivstellensatz and  flat functionals on path $*$-algebras. }

\author{Stanislav Popovych}

\date{}

\maketitle

\begin{abstract}
\medskip{} 
We consider the class of non-commutative $*$-algebras which are path algebras of doubles of quivers with the natural involutions. We study the problem of extending positive truncated functionals on such $*$-algebras.   
An analog of the solution of the truncated Hamburger moment problem \cite{CFHouston} for path $*$-algebras is presented and non-commutative   positivstellensatz is proved.   
We aslo present an analog of the flat extension theorem of Curto and Fialkow for this class of algebras.

\end{abstract}

\section{Introduction.} 

Let $A$ be a $*$-algebra with a generating set $G= \set{x_1, \ldots, x_n, x_1^*, \ldots, x_n^*}$ and assume that there is a linear basis $\B$ consisting of words in the generators  such that $\B \cup \set{0}$ is multiplicatively closed. We will consider $\B$  as linearly ordered with respect to the degree-lexicographic order induced by a linear order on $G$. Then a  {\it moment sequence}  $\alpha = (\alpha_b)_{b\in\B}$ ($\alpha_b \in \Complex$) gives rise to an infinite  matrix $M_\alpha$ with rows and columns indexed by $\B$ such that  
$(M_\alpha)_{b_1,b_2} = \alpha_{b_1b_2^*}$,   which is called the {\it moment matrix} of  $\alpha$. The uppermost left  $n(t)\times n(t)-$corner $M_t$ of $M$ where $n(t)$ is the number of words of length less or equal to $t$  is called the {\it truncated moment matrix} of $\alpha$ of order $t$. The matrix $M_t$ depends only on $\alpha_{2t}$. Thus we can write $M_t = M(\alpha_{2t})$.  We will  denote by $A_t$ the subspace of $A$ generated by the words $b$ of length less or equal to $t$.

When $A$ is the polynomial algebra $\mathbb{R}[x_1, \ldots, x_n]$ and $x_j^*=x_j$  the classical (full) moment problem asks for which $\alpha= (\alpha_b)_{b\in\B}$ there is a representing measure $\mu$, i.e. positive Borel measure on $\mathbb{R}^n$ such that for each monomial $b=b(x_1, \ldots, x_n)$: $$\alpha_b = \int_{\mathbb{R}^n} b(x_1, \ldots, x_n) d \mu(x_1, \ldots, x_n).$$ The question if there is $\mu$ with support contained in set $K\subseteq \mathbb{R}^n$ is called $K-$moment problem.   The classical  Hamburger theorem states that $M$ has a representing measure if and only if  $M$ is positive semidefinite.   However, the truncated moment problem and the truncated $K-$moment problem are more delicate~\cite{curto, curto2, curto3}.  We mention the general solution of the truncated $K-$moment problem in terms of extensions (called truncated version of Riesz-Haviland in~\cite{CFHouston}) and the solution in case of flat $\alpha$ (see~\cite{curto}) because of their relevance to the present paper. 

\newtheorem*{theoremCurt}{Theorem(Curto-Fialkow'2008)}

\begin{theoremCurt}

$\alpha_t$ has $K$-representing measure if and only if  $\alpha_t$ admits $K$-positive extension $\alpha_{t+1}$. 
\end{theoremCurt}
Here $K$-positive means that the  {\it Riesz functional } $L_{\alpha_t}\colon \mathbb{R}[x_1, \ldots, x_n]_t \to \mathbb{R}$ defined as $L_{\alpha_t}(\sum_b \beta_b b) = \sum_b \alpha_b \beta_b$ attains positive values on $p\in \mathbb{R}[x_1, \ldots, x_n]_t$ such that $p(K)\subseteq [0,+\infty)$.

For non-commutative  $*$-algebras $A$ the natural analog of measures are representations on Hilbert spaces. We will say that $\alpha$ (or $M$) is representable if there is a representation $\pi \colon A \to L(H)$ on a pre-Hilbert space $H$ and $\xi \in H$ such that 
$$\alpha_b = \seq{\pi(b)\xi, \xi}.$$    If $J$ is an ideal in $A$ and $\pi(J)=0$ then  we will say that $\alpha$ is $A/J-$ representable.  
It is easy to see that for the non-commutative {\it free} $*$-algebra $\F$  the full moment problem has the same answer as in polynomial algebra case i.e. the full moment problem is solvable iff $M$ is positive semidefinite. An  answer to the truncated moment problem can be given in terms of extensibility (see~\cite{HMPisom}). More precisely, a sequence $\alpha_t$ is representable iff it admits an extension $\alpha_{t+1}$ such that $M_{t+1}$  is positive definite. Moreover, if this is the case,  there is a representation $\pi$ acting on a finite dimensional Hilbert space and a vector $\xi$ such that $\alpha_b = \seq{\pi(b)\xi, \xi}$ for $\abs{b}\le t$.  The same holds for $\F / J$-moment problem for some other ideals $J$, for example, for ideals defining spherical isometries (see~\cite{HMPisom}). The complete analog of this statement is not true for more general ideals.  For instance, there are  ideals $J\subseteq \F$ such that $\F /J$ has only  non-trivial   representations in unbounded operators. For any such  representation $\pi$ the sequence $\alpha_t = (\pi(b))_{\abs{b}\le t}$ clearly has a positive extension but there is 
no representation $\tau$ in {\it bounded} operators with $\tau(J) =\set{0}$ which define $\alpha_t$. To see other possible pitfalls consider the  $*$-algebra $A = \algebra{x, x^* \mid (xx^*)^n}=0$. If $\pi$ is a representation of the free algebra $\F$ and $2t < n$ then the matrix $M_t$ corresponding to  $\alpha_t = (\pi(b))_{\abs{b}\le t}$ in $\F$ is the same as the matrix $M_t$ corresponding to $\alpha_t$ in $A$. Hence there is a  positive matrix $M_t$ for every $t$. Since $A$ has no non-trivial representations the  truncated sequence  $\alpha_t$  is never $A$-representable.

The relation between decompositions of  commutative polynomials into sums of squares and the moment problem is well known. Decompositions of non-commutative polynomials into  sum of hermitian squares in the free $*$-algebras  was obtained in~\cite{Helton} and for some other classes of algebras in  \cite{HMPisom, Schm:AG}.     
In the present paper we study the non-commutative moment problem for  path $*$-algebras $\A_\Gamma$ associated with finite graphs $\Gamma$. Following ideas of~\cite{HMPisom}, we show in Section~3 that a moment sequence $\alpha_t$ is representable if and only if there is a positive semidefinite  extension $\alpha_{t+1}$. Moreover, we show that every hermitian element $f$ which is mapped to a positive semidefinite  operator by every finite dimensional representation of $\A_\Gamma$ is a sum of hermitian squares. Together these facts may be seen as an abstract solution to the non-commutative (full) moment problem for $\A_\Gamma$. 

The other case when the $K$-moment problem for the polynomial algebra has an  especially tractable solution is the case of {\it flat} data $\alpha_t$ which means that $\rank M(\alpha_t) = \rank M(\alpha_{t-1})$. In this case $\alpha_t$ admits a representing measure if and only if  $M_t$ is positive semidefinite~\cite{curto}. Moreover, there is an $r$-atomic representing measure with $r=\rank M(\alpha_t)$.    This is a consequence of the following ``flat extension theorem''. 

\newtheorem*{theoremCF}{Theorem(Curto-Fialkow'1996)}

\begin{theoremCF}   
Every flat $L_{\alpha_t}$ on $\Real[x_1, \ldots, x_n]_{2t}$ can be extended to a flat $L$ on $\Real[x_1, \ldots, x_n]_{2t+2}$.
\end{theoremCF}

In~\cite{laurent, laurent2} M. Laurent has given an algebraic approach to Curto and Fialkow flat extension theorem.   Gr\"obner bases for polynomial algebras appear in several technical constructions in these  papers.     
In Section~\ref{sec:flat} we translate ideas of M. Schweighofer's  unpublished manuscript on Gr\"obner bases  approach to  Curto and Fialkow flat extension theorem  into non-commutative setting. We prove the following analog of Curto and Fialkow flat extension theorem: every flat truncated  functional $L_t$ on $(\A_\Gamma)$ admits a flat extension $L_{t+1}$. In particular, every positive  truncated functional admits an extension to a positive functional on $\A_\Gamma$ and thus via GNS construction defines a finite dimensional representation of $\A_\Gamma$. It is interesting in view of Lance-Tapper conjecture (see Section 2 and~\cite{lance, Tap:JOT_1999}) to obtain an analog of the extension theorem for flat functionals supported on $\A_\Gamma/J$ where $J$ is an ideal generated by paths (i.e. $\A_\Gamma/J$ is a monomial $*$-algebra). In the last section we give some conditions for truncated functionals on free $*$-algebra to have a flat extension.  

 In a forthcoming paper we generalize the extension theorem for a class of $*$-algebras containing all monomial $*$-algebras. In present paper we restrict to the class of path $*$-algebras which makes it possible to minimize the use of  Gr\"obner basis theory and keep formulations of main results  close to classical (commutative) analogues. Note that Gr\"obner basis theory for right ideals is much simpler for  path $*$-algebras than for monomial and commutative polynomial algebras (see~\cite{green}).   

 \section{Definitions.}

Let  $\Gamma = (\Gamma_0, \Gamma_1)$ be a {\it quiver}, that is a directed graph with a finite set of vertices $\Gamma_0$ and  a finite set of arrows $\Gamma_1$. Every arrow $b$ has the unique  origin vertex  $o(b)$ and terminal vertex $t(b)$.  For vertices $e_1, e_2 \in \Gamma_0$ we denote by $\Gamma(e_1, e_2)$ the set of arrows $b$ with $o(b) =e_1$, $t(b) = e_2$. A path $p$ in $\Gamma$ is a finite sequence of arrows (possibly empty) $(v_1, \ldots, v_k)$ such that $o(b_{j+1}) = t(b_j)$ for $j=1, \ldots, k-1$. The number $k$ is called the {\it length} of $p$. The unique empty path at vertex $e$ will be donoted also by $e$. 

Consider {\it the double} $\Gamma^*$ of the graph $\Gamma$ which is the graph with the same set of vertices $\Gamma_0^* = \Gamma_0$ and doubled number of arrows $ \Gamma^*_1 = \Gamma_1\cup  \Gamma_1^*$ where $\Gamma_1^*= \set{ b^* \mid b\in \Gamma_1}$ and if  $b\in \Gamma(e_1, e_2)$ then $b^* \in \Gamma^*(e_2, e_1)$.   

The set consisting of paths $\B$ in a graph $\Gamma$ together with zero element $0$  is a semigroup with multiplication given by concatenation of paths (the product $b_1 b_2$ is $0$ if $t(b_1) \neq o(b_2)$). In particular, for every vertex  $e\in \Gamma_0$ this semigroup contains an idempotent (denoted also by $e$) corresponding to the trivial path at vertex $e$.   The semigroup algebra $\Complex \Gamma$  is called {\it path  algebra } of $\Gamma$ it has a  linear basis consisting of $\B$. The path algebra of $\Gamma^*$ is a $*$-algebra with the involution which maps $b$ to $b^*$.  We will denote this $*$-algebra by $\A_\Gamma$ in sequel.     

Consider  an arbitrary order on the vertices $v_1 < \ldots < v_r$ and arrows $v_r< a_1<a_2 <\ldots$ such that every arrow is greater than every vertex.  Having an order on $\Gamma_0 \cup \Gamma_1$ we can equip the set of paths $\B$ with  the left degree-lexicographic order.   In particular,  for a graph $\Gamma$ with one vertex and $n$ arrows $x_1, \ldots, x_n$ the algebra $\A_\Gamma$ is  the free $*$-algebra $\Complex \seq{x_1, \ldots, x_n, x_1^*, \ldots, x_n^*}$. 
 
Consider the $*$-algebra $$A_{w_1, \ldots, w_m} = \Complex \seq{x_1, \ldots, x_n, x_1^*, \ldots, x_n^* \mid w_1 =0, \ldots, w_m =0,  w_1^* =0, \ldots, w_m^* =0}$$ where $w_j$ are words in $x_j, x_j^*$.  We can always assume that the set $\{ w_1, \ldots, w_m,$ $w_1^*, \ldots, w_m^*\}$ is reduced, i.e.  no word is a subword  of some other word. Such algebras $A_{w_1, \ldots, w_m}$ constitute a subclass of {\it monomial} algebras. A word $w$ is called {\it unshrinkable} if $w$ can not be presented as $u d d^*$ or $d d^* u$ for some word $u$ and non-empty word $d$. We call $A_{w_1, \ldots, w_m}$ a Lance-Tapper $*$-algebra if every $w_j$ is unshrinkable. It was conjectured by Lance and Tapper that a $*$-algebra $A_{w_1, \ldots, w_m}$ has a faithful representation on a Hilbert space if and  only if it is a Lance-Tapper $*$-algebra (see~\cite{lance, Tap:JOT_1999, Pop:monom}). 
Moreover, they conjectured that every Lance-Tapper $*$-algebra has a  separating family of finite dimensional representation. The following simple lemma shows that conjecture is true for trivial monomial ideals.

\begin{lemma}\label{faithrep} For any quiver $\Gamma$ the path $*$-algebra $\A_\Gamma$ has a separating family of finite dimensional representations. 
\end{lemma}
\begin{proof}
For a path $*$-algebra $\A_\Gamma$ there is a $*$-isomorphism $\phi$ of $\A_\Gamma$ with a subalgebra in $M_n(\F)=\F \otimes M_n$ were $n = \abs{\Gamma_0}$ and $\F$ is a free $*$-algebra with free generators $x_b$ corresponding to each $b\in \Gamma_1$. Let $e_1, \ldots, e_n$ be enumeration of all vertices of $\Gamma$. Under the isomorphism $\phi$ arrow $b\in\Gamma(e_i,e_j)$ maps to $x_b\otimes E_{ij}$ and idempotent $e_i$ maps to $E_{ii}$. Since the free algebra $\F$ has a separating family of finite dimensional representation (see~\cite{Avitzour, Pop:monom}) the same is true for any path $*$-algebra. 
\end{proof}

There is a well known bijective correspondence between representation of $\A_\Gamma$ and representations of $\Gamma$ in Hilbert spaces. Recall that a representation $\Pi$ of $\Gamma$ in a category of Hilbert spaces  is a mapping which maps each vertex $e\in \Gamma_0$ into a Hilbert space $H_e$ and each arrow $b\in \Gamma(e_1, e_2)$ into a linear operator $\Pi(b)\colon H_{e_1} \to H_{e_2}$.  To a representation $\Pi$ of $\Gamma$  there corresponds a representation $\pi$  of $\A_\Gamma$ in the Hilbert space $H = \oplus_{i=1}^n H_{e_i}$ such that for $b\in \Gamma(e_i,e_j)$ the operator  $\pi(b) \colon H \to H$ is given by the block--matrix $\Pi(b)\otimes E_{ij}$ and $\pi(b^*)$ is defined as the adjoint operator to $\pi(b)$. Conversely, given a representation $\pi$ of $\A_\Gamma$ in a Hilbert space  we put $H_e = \pi(e)(H)$ (here $e$ is the trivial path at vertex $e$). Since for $b\in \Gamma(e_1, e_2)$, $e_1 b e_2 = b$ one can check that the operator $\pi(b)$ maps subspace $H_{e_2}$ into $H_{e_1}$. Thus mapping $e\mapsto H_e$, $b \mapsto  \pi(b)|_{H_{e_2}}$ defines a representation of $\A_\Gamma$. 

\section{Sums of squares decomposition. }

Enumerating the vertices $\Gamma_1$ as $b_1, \ldots, b_n$ we get a one-to-one correspondence between the class of finite dimensional representations $\Rep_{f.d.} \A_\Gamma$ of $\A_\Gamma$ and the class $Z$ of  $n-$tuples of operators $(X_1, \ldots, X_n)$ acting on a finite dimensional Hilbert spaces such that 
\begin{gather} \label{pathrel1}
X_i X_j =0 \text{  if } b_i b_j =0,\\  
X_i^2 =X_i \text{ if } b_i =b_i^2,\\
X_i X_j = X_i \text{ if } b_i b_j =b_i,\\
X_i X_j = X_j \text{ if } b_i b_j =b_j.  \label{pathrel2}
\end{gather}  Clearly, $\A_\Gamma$ is isomorphic to the quotient of the free $*$-algebra $\F$ with the generating set $X = \Gamma_1\cup\Gamma_1^*$ by the $*$-ideal generated by  relations~\eqref{pathrel1}--\eqref{pathrel2}.  
Clearly,  for any unitary operator $U \colon H \to K$ between two Hilbert spaces $n-$tuple $(U A_1 U^*, \ldots, U A_n U^*) \in Z$. Hence $Z$ is closed under joint unitary transformation. 
The following lemma is a direct consequence of~\cite{HMPisom} and Lemma~\ref{faithrep}.  Recall that $(\A_\Gamma)_d$ denote the subspace generated by paths of length no greater than $d$. 
\begin{lemma}\label{closedcone} For any $d\ge 1$ the cone  $\mathcal{C}_{2d}(\A_\Gamma)=\co \set{ff^* \mid f\in (\A_\Gamma)_d}$ is closed in $(\A_\Gamma)_{2d}$. Here $\co$ denote convex hull.   
\end{lemma}    
\begin{proof}
The ideal $I(Z)$ in the free algebra $\F$ consisting of $p(x_1, \ldots, x_n, x_1^*, \ldots, x_n^*)$ such that $p(X_1, \ldots, X_n, X_1^*, \ldots, X_n^*) = 0$ for all $(X_1, \ldots, X_n)\in Z$ coincides with the kernel of the canonical surjection $\psi \colon \F \to \A_\Gamma$  by Lemma~\ref{faithrep}. 
The image  $\psi(\co \set{ff^* \mid f\in \F_d})$,  denoted by $\mathcal{C}_{2d}(Z)$ in~\cite{HMPisom}, is closed in $\mathcal{P}_{2d} = \psi(\F_{2d})$ by~\cite[Lemma~3.2]{HMPisom}.    
Identifying $\F / I(Z)$ with $\A_\Gamma$ via $\psi$ we have that $\mathcal{P}_{2d} = (\A_\Gamma)_{2d}$ and  $ \mathcal{C}_{2d}(Z) = \mathcal{C}_{2d}(\A_\Gamma)$. Hence, the lemma follows. 

\end{proof}

\begin{lemma}\label{compr}
 
Let $H$ be a pre-Hilbert space, $d\ge 1$ and $H_0 \subseteq H_1 \subseteq \ldots \subseteq H_{d}= H$ be subspaces with $H_0$ finite dimensional. Assume that for  $0\le t\le d-1$ we are given a linear map $\pi_t \colon (\A_\Gamma)_{t} \to L(H_{d-t}, H)$ such that 
\begin{enumerate}
\item for every $f\in (\A_\Gamma)_{t}$ and $0\le s\le d-t$,  $\pi_t(f)(H_{d-t-s}) \subseteq H_{d-s}$.
\item for every $0\le r \le t$ and  $g\in (\A_\Gamma)_{t-r}$ we have 
$\pi_{t-r}(g)|_{H_{d-t}}= \pi_{t}(g)$. Thus we can omit subscript $t$ in the notation $\pi_t(f)$. 
\item  for every $f_1\in (\A_\Gamma)_{t_1}$, $f_2\in (\A_\Gamma)_{t_2}$ with $t= t_1+t_2\le d-1$ we have  $\pi(f_1 f_2)|_{H_{d-t}} = \pi(f_1) \pi(f_2)|_{H_{d-t}}$
\item For $f\in (\A_\Gamma)_{t}$ and $u, v \in H_{d-t}$,  $$ \seq{\pi(f)u,v} = \seq{u,\pi(f^*) v}.$$
\end{enumerate}
Then  there is a finite dimensional subspace $H'\subseteq H$ with  $\dim H' \le \dim H_0 \dim (\A_\Gamma)_d$,  $H_0\subseteq H'$ and a representation $\tau$ of $\A_\Gamma$ on $H'$ such that $$ \pi(a)|_{H_0} = \tau(a)|_{H_0} \text{ for all } a \in (\A_\Gamma)_{d-1} .$$

\end{lemma} 
\begin{proof} 
Put $H' =\set{ \pi(f) H_0 \mid f\in (\A_\Gamma)_{d}}$ and $K = \set{ \pi(f) H_0 \mid f\in (\A_\Gamma)_{d-1}}$ and $V= \set{ \pi(f) H_0 \mid f\in (\A_\Gamma)_{d-2}} $ which are finite dimensional subspaces of $H$. We will define  a representation $\Pi$ of $\Gamma$. For $e\in \Gamma_0$ $\pi(e)$ is a projection defined on $H$. Put $H_e = \pi(e) H'$ and $K_e = \pi(e) K$. Let $K_e^\bot$ be the orthogonal complement of $K_e$ in $H_e$, i.e. $H_e = K_e \oplus K_e^\bot$. 

For $b\in \Gamma_1(e_1, e_2)$ put $\Pi(b)|_{K_{e_2}} = \pi(b)|_{K_{e_2}}$ and $\Pi(b)|_{K_{e_2}^\bot} =0$. Since $\pi(e_1) \pi(b) = \pi(b)$ on $H_{d-1}$ we have $\Pi(b) \colon H_{e_2} \to H_{e_1}$ and $K_{e_2}^\bot \subseteq \ker \Pi(b)$. Hence for the adjoint operator $\Pi(b)^*$ we have  $$\Pi(b)^* \colon H_{e_1} \to H_{e_2} \text{ and } \Ran\Pi(b)^* \subseteq K_{e_2}.$$ 
From this follows that 
\begin{equation}\label{grapheq}
\pi(b^*)|_{K_{e_1}\cap V} = \Pi(b)^*|_{K_{e_1}\cap V}.
\end{equation} 
Indeed, for any $v \in K_{e_1}\cap V$ and $w\in K_{e_2}$ we have 
$$ \seq{w, \Pi(b)^* v } = \seq{\Pi(b)w, v} = \seq{\pi(b)w, v} =\seq{w, \pi(b^*) v}.$$  
Since $\Pi(b)^* v \in K_{e_2} $, $\pi(b^*) v \in K_{e_2}$ and $w\in K_{e_2}$ is arbitrary we get $\Pi(b)^* v =  \pi(b^*) v$.  Let $\tau$ be the representation of $\A_\Gamma$ corresponding to $\Pi$. Then $\tau$ is defined on the Hilbert space  $\oplus_e H_e$ which we will identify with $H'$. By~\eqref{grapheq} and induction in the length of path $b\in \B$ we get that $\tau(b)|_{H_0} =\pi(b)|_{H_0}$ for all paths $b$ with $\abs{b}\le d-1$.  
\end{proof}
\begin{remark}
If $\pi$ is a (possibly unbounded) representation on a pre-Hilbert space $H$ and $H_0\subseteq H$ is a finite dimensional subspace then  with $H_t = \pi((\A_\Gamma)_t) H_0$ we have for $f\in (\A_\Gamma)_t$,  $\pi(f)(H_{d-t})\subseteq H_d$ and conditions (1)-(4) are satisfied. Hence there is finite dimensional representation 
$\tau$ such that $\pi(f)|_{H_0} = \tau(f)|_{H_0}$ for $f\in (\A_\Gamma)_d$.     
\end{remark}

The following corollary is an analog of the solution of the  truncated  Hamburger moment problem from~\cite{CFHouston}. 
\begin{corollary} 
A truncated linear functional $L_{d}$ on path $*$-algebra $(\A_\Gamma)_{2d}$ admits an  extension to a positive functional $L$ on $\A_\Gamma$ if and only if it admits and extension to positive functional $L_{d+1} \colon (\A_\Gamma)_{2d+2} \to \Complex$.   
\end{corollary} 
\begin{proof} Let $K =\set{f \in (\A_\Gamma)_{2d+2} \mid L_{d+1}(ff^*) =0}$. Since $L_{d+1}$ is positive $\abs{L_{d+1}(fg^*)} \le 
{L_{d+1}(ff^*)}^{1/2} {L_{d+1}(gg^*)}^{1/2}$ and $K$ is the null space of the sesquilinear for $\seq{f,g} = L_{d+1}(fg^*)$ on the space $(\A_\Gamma)_{2d+2}$. Then $\seq{\cdot, \cdot}$ induces  an inner product on quotient space $H= (\A_\Gamma)_{d+1} /K$. Let $H_t = (\A_\Gamma)_{t} +K \subseteq H$. Let  $H_0$ be one dimensional subspace generated by  $\xi = 1+K$. Given $f\in (\A_\Gamma)_{t}$ define $\pi_t(f)$ as the restriction of right multiplication by $f$ operator on $H_{d+1-t}$.   It is routine to check conditions (1)-(4) of  Lemma~\ref{compr}.  Thus there is a finite dimensional representation $\tau$ such that for all $f,g\in (\A_\Gamma)_{d}$ we have  $L_{d}(fg^*) = \seq{ \pi_d(f)\xi, \pi_d(g)\xi  } = \seq{\tau(f)(\xi), \tau(g)\xi} = \seq{\tau(fg^*)(\xi), \xi}$. Which proves that $L$ has the positive extension $\seq{\tau(\cdot)(\xi), \xi}$.

\end{proof}

Let $\Sym_d(\Gamma)$ be the set of hermitian elements $f$ of $\A_\Gamma$ with $\deg(f)\le d$. 
By  straightforward reformulation of  \cite[Lemma~3.1]{HMPisom} combined with  Lemma~\ref{compr} we get the following lemma.  
\begin{lemma} \label{biort}
For every $d\ge 0$ there is a basis $\beta_1$, $\ldots$, $\beta_k$ of $\Sym_d(\Gamma)$, a sequence  $\pi_1$, $\ldots$, $\pi_k$ of representations of  $\A_\Gamma$ on a finite dimensional Hilbert space $H$ and a vector $\xi \in H$ such that for all $i,j$:  \begin{equation} 
\seq{\pi_i(\beta_j)\xi, \xi} = \delta_{ij},
\end{equation} 
where $\delta_{ij}$ is Kronecker's symbol. 
\end{lemma}

Following~\cite{HMPisom} we can introduce a norm $\norm{\cdot}$ on $\Sym_d(\Gamma)$ by the formula 
\begin{equation}\label{HMPnorm}
\norm{f} = \sum_{i=1}^k \abs{ \seq{\pi_i(f) \xi, \xi} }
\end{equation}
for $f\in \Sym_{d}(\Gamma)$. For even $d$ and any $h\in \C_{d}(\Gamma)$, $\norm{h}$ can be expressed as a value of the linear functional $N(h)=\sum_{i=1}^k \seq{\pi_i(h) \xi, \xi}$, i.e. $\norm{h}= N(h)$.  

The following theorem shows that a hermitian element of $\A_\Gamma$ which is positive semidefinite in every finite dimensional representation is a sum of hermitian squares. 
  \begin{theorem}\label{sosrep}
  Let a hermitian $q\in (\A_\Gamma)_{d-1}$ be such that for any representation $\pi$ of $\A_\Gamma$ with $\dim \pi\le \dim (\A_\Gamma)_d$, $\pi(q)$ is positive semidefinite then $q\in \C_{2d}$.   
  \end{theorem}
  \begin{proof}
  Assume that $q\not\in \C_{2d}$. Since $\C_{2d}$ is closed in $\Sym_d(\Gamma)$   Minkowski's separation theorem implies that there is a linear functional $L_0 \colon (\A_\Gamma)_{2d} \to \Complex$ such that $L_0(q) <0 \le \min_{c \in \C_{2d}}{  L_0(c)}$. Take $\epsilon>0$ then for $L =L_0+\epsilon N$ we have $L(ff^*)>0$ for $f\in (\A_\Gamma)_d \setminus \set{0}$. Hence $\seq{f,g} = L(fg^*)$ defines a scalar product on $H=(\A_\Gamma)_d$. 
  
Let $H_0 =\Complex 1$ and  $H_t = (\A_\Gamma)_t$. For $g\in (\A_\Gamma)_t$ let $R_g:H_{d-t} \to H_d$ denote the operator of multiplication by $g$ from the right. By Lemma~\ref{compr} there is a finite dimensional representation in Hilbert space $H'$ such that $H_0 \subseteq H'$ and $\dim H' \le \dim (\A_\Gamma)_d$ such that  $\tau(g) 1 = R_g 1 = g$ for all $g\in (\A_\Gamma)_{d-1}$. Hence $L(ff^*)= \seq{R_f 1,1} =\seq{\tau(f)1,1} \ge 0$.        
  
  \end{proof} 

\section{Flat functionals.}\label{sec:flat}
In this section we study  which truncated positive semidefinite functional on a path $*$-algebra can be extended to a positive semidefinite functional on the whole algebra with the moment matrix of finite rank. 
We need some standard definitions from  Gr\"obner basis theory. 
We will consider associative algebras over field of complex numbers. 
For an algebra $A$ to have Gr\"obner basis theory, $A$ must have a  multiplicative linear basis $\B$ (i.e. for every $b_1$, $b_2 \in \B$, $b_1 b_2 \in \B$ or $b_1 b_2 =0$) with an admissible order on $\B$ (see~\cite{green}).    
An order $>$ is called {\it admissible} if 
\begin{enumerate}
\item[A0.] $>$ is well-order  on $\B$.  
\item[A1.] For all $b_1, b_2, b_3 \in \B$, if $b_1 > b_2$ then $b_1 b_3 > b_2 b_3$ if both $b_1 b_3$ and $b_2 b_3$ are nonzero.  
\item[A2.] For all $b_1, b_2, b_3 \in \B$, if $b_1 > b_2$ then $b_3 b_1 > b_3 b_2$ if both $b_3 b_1$ and $b_3 b_2$ are nonzero.
\item[A3.]  For all $b_1, b_2, b_3, b_4 \in \B$, if $b_1 = b_2 b_3 b_4$ then $b_1 \ge  b_3$. 
\end{enumerate}
Following~\cite{green} we say that $A$ has an ordered multiplicative basis $(\B, >)$ if $\B$ is a multiplicative basis and $>$ is an admissible order on $\B$. 
It was shown in~\cite{green} that every algebra with  ordered multiplicative basis is a quotient of a path algebra by  $2-$nomial  ideal. By {\it $2-$nomial} ideal we mean an ideal generated by some  elements of the form $p$ or $p-q$ where $p$, $q\in \B$. An algebra $\A$ which is a quotient of a path algebra by the ideal generated by a set of paths is called {\it monomial}.

 Let $x = \sum_{i=1}^r \alpha_i b_i$ were $\alpha_i \in \Complex^*$  and $b_i$ are distinct elements of $\B$.  The tip of $x$, denoted by $\Tip(x)$, is the largest element in $\set{b_1, \ldots, b_r}$, i.e. $\Tip(x) = b_j$ for $b_j$ such that $b_j\ge b_i$ for all $i=1, \ldots, r$.   
 
 For  $k\ge 1$ put $S_k = \set{ a\in \B \mid \abs{a} \le k}$.  Denote by $V_k$ the linear span of $S_k$ and by  $V_k^*$ the dual vector space of $V_k$. Note that  $V_k$ is the same thing as $A_k$ in the preceeding sections. We change the notations to avoid ambiguous notation $A_k^*$.  
 
\begin{definition} 
Given $L_k \in V_{2k}^*$ define $B_{L_k} \colon V_k \times V_k \to \Complex$ by 
$B_{L_k}(p, q) = L_k(p q^*)$. Functional $L_k$ with $k>1$ (and sesquilinear form $B_{L_k}$) will be called a flat truncated functional on $A$ (resp. flat sesquilinear form on $A$)  if $L_k$ is hermitian (i.e.  $L_k(a^*)= \overline{L_k(a)}$ for all $a\in A$) and $\rank B_{L_k} = \rank B_{L_{k-1}}$ with $L_{k-1}$ being  the restriction of $L_k$ to $V_{k-1}$. 
\end{definition}
Note that the sesquilinear form $B_{L_k}$ is given by the moment matrix $M_{\alpha}$ (were $\alpha_b =L_k(b)$) in the basis $S_k$ of vector space $V_k$.   

Clearly $\ker B_{L_k} \cap V_{k-1} \subseteq \ker B_{L_{k-1}}$ and we have linear maps $\pi$ and $i$ such that 
 
\begin{equation}\label{diag}
\xymatrix{ V_{k-1} / \ker B_{L_{k-1}}   &
V_{k-1}/(\ker B_{L_{k}} \cap V_{k-1}) \ar@{->>}[l]_{\ \pi   \quad } \ar@{^{(}->}[r]^{ \quad \quad i  } & V_{k} / \ker B_{L_{k}}
  }
\end{equation}

In the above diagram dimension is weakly decreasing from left to right. Thus $L_k$ is flat implies that 
$\pi$ and $i$ are isomorphisms. Hence $L_k$ is flat if and only if  
\begin{gather} 
\ker B_{L_k} \cap V_{k-1} = \ker B_{L_{k-1}},\label{flat1} \\
V_{k} = V_{k-1} +\ker B_{L_k}.\label{flat2}
\end{gather}   

We have the following analog of {\it recursively generated}  property of $\ker B_{L_k}$ studied  for commutative  polynomial algebras in~\cite{curto3}. 
\begin{lemma} \label{trunk}
Let $A=\A_\Gamma$ be a path $*$-algebra. 
If $L_k$ is flat then for every $p\in \ker B_{L_k}$ and $w\in \B$ such that 
$\Tip(p)w \neq 0$ and $p w \in V_k$ we have $p w \in \ker B_{L_k}$. 
\end{lemma}
\begin{proof} 
We can assume that $w$ belongs to the  set $S_1$ (otherwise decompose $w= w_1 \ldots w_s$ and use induction on $s$).  
Given $v \in V_k$  we can decompose $v = h + g$ where $h\in V_{k-1}$, $g\in \ker B_{L_k}$. 
Hence $$L_k(p w v^*) = L_k(p (w  h^*)) + \overline{L_k(g (p w)^*)} = L_k(p (w  h^*))= 0.$$ 
The second equality follows from  $g\in \ker B_{L_k}$ and $p w \in V_k$  and the last one follows from $p \in \ker B_{L_k}$ and $w h^* \in V_k$. 

\end{proof}

We will show that any flat truncated functional $L_k$  on a path $*-$algebra $A$ 
can be extended to a hermitian functional on $A$. For this we need basic facts about right  Gr\"obner bases for right ideals.   
There is a general theory of right Gr\"obner basis for right  modules with coherent bases  over algebras  with ordered multiplicative basis~\cite{green}. Recall that a linear basis   $\mathcal{M}$   of  a right module $M$ over an algebra $A$ with ordered multiplicative basis $(\B, >)$ is {\it coherent} if for all $m\in \mathcal{M}$ and every $b\in \B$ either $m b \in \mathcal{M}$ or $m b = 0$. A well-order $\succ$ on $\mathcal{M}$ is a {\it right admissible order on } $\mathcal{M}$ if 
\begin{enumerate}
\item  For all $m_1, m_2 \in \M$ and $b\in \B$, if $m_1 \succ m_2$ then 
$m_1 b \succ m_2 b$ if both $m_1 b$ and $m_2 b$ are nonzero.  
\item  For all $m \in \M$ and $b_1, b_2\in \B$, if $b_1 > b_2$ then 
$m_1 b \succ m_2 b$ if both $m b_1$ and $m b_2$ are nonzero.
\end{enumerate}

If $x \in \M \setminus\set{0}$ then $x= \sum_{j=1}^s \alpha_j m_j$ were  $\alpha_j \in \Complex^*$ and $m_1, \ldots, m_s$ are distinct elements of $\M$. The tip of $x$ is $\Tip(x) =m_i$ such that  $m_j \succeq m_i$ for all $i=1, \ldots, s$.   

Recall the algorithm from~\cite{green} of constructing a right Gr\"obner basis for a submodule $M$   with a coherent ordered basis $(\mathcal{M}, \succ)$ of a right projective module over an path $*$-algebra $A$ with an  ordered multiplicative basis $(\B, >)$. This algorithm will be needed in the proof of Theorem~\ref{gbt}.  Assume that $M$ is generated by a finite set of elements $H= \set{h_1, \ldots, h_n}$. For a subset $X \in M$ denote $\Tip(x) =\set{\Tip(x) \mid x\in X}$ and $\NonTip(X) = \M \setminus \Tip(X)$.  If $\mathcal{N}$ is a right submodule of $\M$ then a subset $\G\subset \N$ is called a {\it right Gr\"obner basis} of $\N$ if $\Tip(\G)$ generate $\N$ as a submodule. Consider the following algorithm of transforming the subset $H$. 

\begin{enumerate} 
\item Remove $0$ from $H$.   
\item Put $\mathcal{T}_H = \set{\Tip(h) \mid \text{ for all } h' \in H \setminus \set{h}, \   \Tip(h') \text{ does not left divide } \Tip(h) }$. 
\item For every $t \in \mathcal{T}_H$ choose $h\in H$ such that $\Tip(h) =t$ and renumber so that these elements are $h_1, \ldots, h_s$. If $s=n$ we are done. Otherwise denote by $Q^\dagger$ the right submodule generated by $h_1, \ldots, h_s$. 
\item For every $i=s+1, \ldots, n$ decompose $h_i = h_i^\dagger + \Norm(h_i)$ using  
$M = Q^\dagger\oplus \Span (\NonTip(Q^\dagger))$, i.e. using $h_1, \ldots, h_s$ reduce $h_i$. After finite steps of reductions we get $\Norm(h_i)$.
\item Put $H= \set{ h_1, \ldots, h_s, \Norm(h_{s+1}), \ldots, \Norm(h_{n})}$. 
\end{enumerate} 

Here reduction means the following. If $h_i = \sum \alpha_i m_i$, $\alpha_i \in \mathbb{C}$, $m_i \in \mathcal{M}$ and for some $i$ and some $k\in \set{1,\ldots,s}$, $\Tip(h_k)$ left divides $m_i$ ( i.e. there is $b\in \B$, $\alpha\in \Complex$ such that $m_i = \alpha \Tip(h_k) b$) then we can make a reduction which means that  we replace $m_i$ in the decomposition of $h_i$ above  by the element  $m_i - \alpha h_k b$.  The total reduction of $h_i$ by $H$ is a sequence of reductions by $H$ that can not be further reduced by $H$. It was proved in~\cite[Proposition 4.2]{green} that  any finitely generated submodule $\N$ of a projective module over a path algebra has a right Gr\"obner basis which can be computed by the above algorithm.

An algebra $A$ with an ordered multiplicative basis $(\B, >)$ is a right projective  module over $A$ with respect to the  multiplication from the right and putting $\M =\B$   we get a coherent module which has an admissible order $\succ$ equal to $>$.  We will denote this module by $A_r$. In particular, every  right ideal $J$ is a submodule of $A_r$ and hence has a right Gr\"obner basis.

\begin{theorem} \label{gbt}
Let $L_k \in V_{2k}^*$ be a flat functional on a path $*$-algebra $A= \A_\Gamma$. Let $J$ be the right ideal generated 
by $\ker B_{L_k}$ in $A$.  Then there is a right Gr\"obner basis $\mathcal{G}$ such that 
 $\mathcal{G} \subseteq \ker B_{L_k}$. 
\end{theorem}
\begin{proof} 
Take any generating set $\set{h_1, \ldots, h_n}$ of $\ker B_{L_k}$. 
If $\Tip(h_i)$ left divides $\Tip(h_j)$, i.e. $\Tip(h_j) = \Tip(h_i) b$ where $b\in \B$ then $h_j -h_i b \in \ker B_{L_k}$ since $h_i b \in \ker B_{L_k}$ by Lemma~\ref{trunk}. This proves that step (4) of the algorithm of computing rigth Gr\"obner basis transforms elements of $\ker B_{L_k}$ into elements of $\ker B_{L_k}$. All other steps of the algorithm,  clearly, produce no new elements. Thus obtained right Gr\"obner bases will be contained in $\ker B_{L_k}$. 

\end{proof}

\begin{theorem} \label{thmExt}
Any flat truncated functional  $L_k$  on a path $*$-algebras $A=\A_\Gamma$ has an extension 
to a linear functional on $A$ such that \begin{equation}\label{rk}
\rank B_{L_k} = \rank B_{L}.
\end{equation}
 An extension which satisfies~\eqref{rk} is unique and positive semidefinite if $L_k$ is such.  
\end{theorem}
\begin{proof}
Let $J$ be the rigth ideal generated by $\ker B_{L_k}$ in $A$.  Let us denote by $\widetilde{B_{L_k}}$ the induced sesquilinear form on $V_k / \ker B_{L_k}$. 
The natural map $\phi \colon V_k / \ker B_{L_k} \to A/ J$ is linear. The map $\phi$ is onto since modulo $\ker B_{L_k}$ each element of $V_{k}$ is equivalent to en element of $V_{k-1}$ and hence, by induction,  modulo right ideal $J$ each element of $A$ is equivalent  to some element from $V_{k-1}$. The same inductive argument also shows that for every $b\in S_{2k}$ there are $v\in V_{k-1}$, $g_t \in \ker B_{L_k}$ and $b_t \in S_k$ ($t =1, \ldots, l$) such that 
\begin{equation}\label{recform}
b= v+ g_1 b_1 +\ldots + g_l b_l.
\end{equation} 
Hence 
\begin{equation}\label{resp}
L_k(b)= L_k(v)+ B_{L_k}(g_1, b_1^*) +\ldots + B_{L_k}(g_l, b_l^*) = L_k(v).
\end{equation}
 
By Theorem~\ref{gbt} we can find a   a Gr\"obner basis $\G$ for the right ideal $J$  such that $\G \subseteq \ker B_{L_k}$. If $x \in V_k$ such that $\phi(x) = 0$ (i.e. $x \in J$ ) then since $J$ is a right Gr\"obner basis \begin{equation}\label{form} 
x= \sum_{j=1}^s \alpha_j g_j b_j
\end{equation} where $\alpha_j \in\Complex^*$, $g_j \in \G$ and $b_j \in \B$. Moreover, we can  assume without loss of generality that $\Tip(g_j) b_j \neq 0$ and $g_j b_j \in V_k$. Indeed, by the definition of a right Gr\"obner basis, there is $g \in \G$ such that $\Tip(x) = \alpha \Tip(g) b$ for some $b \in \B$ and $\alpha\in \Complex^*$. Thus  $\Tip(g) b \neq 0$, the element $y= x - \alpha g b \in J$ and  $\Tip(y)$ strictly less than $\Tip(x)$. Hence, by induction,  we get that $x$ has required decomposition~\eqref{form} with  $\Tip(g_j) b_j \neq 0$ for all $j$. 
By Lemma~\ref{trunk} we get that $g_j b_j\in \ker B_{L_k}$ for all $j=1, \ldots, s$ and consequently  $x \in \ker B_{L_k}$. Thus $\phi$ is injective. 
 By \eqref{recform} we also have 
 \begin{equation}\label{extend}
 L_k( p ) = L_k(\phi^{-1}(p+J)) \text{ for all } p \in V_{2k}. 
\end{equation}
Thus putting $L(p) = {L_k}(\phi^{-1}(p+\J))$ for $p\in A$ we get an extension of  $L_k$.  

Let us show that $\ker B_{L} = J$. For every $p\in J$ and  $q\in A$ we have $pq \in J$, thus $B_{L}(p,q)=0$. This implies that  $J\subseteq \ker B_{L}$. For the converse inclusion consider $q\in \ker B_L$ and take $p\in V_k$ such that $p-q \in J$. Since $J\subseteq \ker B_{L}$ we have $p\in \ker B_{L} \cap V_k$. Hence $p \in \ker B_{L_k}\subseteq J$ and we conclude that $\ker B_L =J$. 

For any $p_1, p_2 \in A$ take $q_1, q_2 \in V_k$ such that the elements $d_j= p_j-q_j$ belong to $V_k$. Then $L(d_1 p_2^*) = 0$, which implies $L(p_1 p_2^*) = L(q_1 p_2^*)$. 
Analoguosly $L(q_1 d_2^*) =0 $, which implies $L(q_1 p_2^*) = L(q_1 q_2^*)$ and we conclude that $L(p_1  p_2^*) = L_k(q_1 q_2^*)$. From this follows that $L$ is hermitian and positive definite if $L_k$ is such.    

Let us prove uniqueness. Let $L' \in A^*$ be a linear extension of $L_k$ to $A$ such that $\rank B_{L'} = \rank{B_{L_k}}$. The subspace $\ker B_{L'} \cap V_k =\set{v\in V_k \mid L'(v p^* ) = 0 \text{ for all } p\in A}$ is clearly contained in $\ker B_{L'_k}$ which is equal to $\ker B_{L_k}$ since $L$ is an extension of $L_k$. Thus we have the following diagram  
\begin{equation}\label{diag}
\xymatrix{ A / \ker B_{L'}   &
V_{k}/ \ker B_{L'} \cap V_k \ar@{->>}[r]^{  \pi  \quad \quad  } \ar@{_{(}->}[l]_{  i \quad  } & V_{k} / \ker B_{L_{k}} \cong A / \ker B_L 
  }
\end{equation}
The dimensions in the above diagram is decreasing from left to right and we get that $\ker B_{L'} \cap V_k = \ker B_{L_k}$ and $i$ is  surjective. 
From this follows that for any $p\in A$ there exists $q\in V_k$ such that $p-q \in \ker B_{L'}$. Hence $L'(p) =L'(q) = L_k(q)$. 
\end{proof}

\begin{lemma} \label{lemF}
Let $L_k \in V_{2k}^*$ be hermitian. Decompose   $B_{L_{k}}$  into the block-matrix 
$$\left(\begin{array}{cc}
A & C  \\
{C}^* & B \\
\end{array}\right)
$$ with respect to  the 
decomposition $V_k = V_{k-1}\oplus \Span(S_k\setminus S_{k-1})$. 
Then $L_k$ is flat if and only if $$\Ran(C)\subseteq \Ran(A) \text{  and }    B= C^*(A|_{\Ran(A)})^{-1} C.$$ 
If $B_{L_{k-1}}$ is positive semidefinite then so is $B_{L_{k}}$. 
\end{lemma} 
\begin{proof} Assume first that $L_k$ is flat. 
By~\eqref{flat2} for every  $b\in S_k\setminus S_{k-1}$ there is $v\in V_{k-1}$ such that $b-v \in \ker B_{L_{k}}$ and hence  
\begin{gather*} 
A v = C b, \\
C^* v = B b.
\end{gather*}
 Since the set of $b$ generates $\Span(S_k\setminus S_{k-1})$ which is a domain of $C$  the first equation above means that  $\Ran(C) \subseteq \Ran(A)$. If $\Ran(C) \subseteq \Ran(A)$ then the operator $(A|_{\Ran(A)})^{-1} C$ is well defined and $v= (A|_{\Ran(A)})^{-1} C b$. Hence $B b= C^*(A|_{\Ran(A)})^{-1} C b$. Thus  $B= C^*(A|_{\Ran(A)})^{-1} C$. 

Assume now that $\Ran(C)\subseteq \Ran(A)$ and   $B= C^*(A|_{\Ran(A)})^{-1} C$. 
 For every $b\in S_k\setminus S_{k-1}$ we can put  $v= (A|_{\Ran(A)})^{-1} C b$.  It follows that $V_k = \ker B_{L_k} +V_{k-1}$.  We need to show that $ \ker B_{L_{k-1}}  \subseteq \ker B_{L_{k}}$. 
 
Let $\seq{\cdot, \cdot}$ denote the inner product on $A$ such that  $\B$ is an orthonormal basis.  In particular $A_{b_1, b_2} = \seq{b_1, A b_2}$.
 If $v \in \ker B_{L_{k-1}}$ then $A v = 0$ hence $\seq{v, A v'} = 0$ for all $v'\in V_{k-1}$. Hence 
$v$ is orthogonal to $\Ran(A)$ and consequently to $\Ran(C)$. Hence $\seq{v , C w} =0$ for all $w\in \Span(S_{k}\setminus S_{k-1})$. Thus  $B_{L_k}(v,w) = \seq{v,C w} =0$ and $\ker B_{L_{k-1}}  \subseteq \ker B_{L_{k}}$. 
Form this and    $V_k = \ker B_{L_k} +V_{k-1}$ we get that $\dim \ker B_{L_k} =\dim \ker A + \dim \Span(S_{k}\setminus S_{k-1})$. From this follows that nonzero eigenvalues of $B_{L_k}$ are exactly  nonzero eigenvalues of the matrix $A$.

\end{proof}

 For the class of free $*$-algebras some flat functionals can be obtained as extensions of tip--maximal functionals defined below. It is easier to construct tip--maximal functionals which will be used in examples. 
 
\begin{definition}
Let us call a hermitian functional $L_k$ {\it tip--maximal} if  there is a  generating set $p_1, \ldots, p_r$ for $\ker B_{L_k}$. 
Such that projections the $f_k$ of $p_k$ on  $\Span(S_k\setminus S_{k-1})$ parallel to $\Span(S_{k-1})$  is linearly independent. 
 In particular, $L_k$ is tip--maximal if $B_{L_k}$  is non-degenerate.  
\end{definition}

\begin{lemma}\label{lemRec}
Every tip--maximal $L_{k-1}$ on a free $*$-algebra can be extended to a flat $L_{k}$.  Moreover, if $L_{k-1}$ is positive semidefinite then $L_k$ is such. 
\end{lemma} 
\begin{proof} 
We will construct matrix of $L_{k}$ in the form $$\left(\begin{array}{cc}
A & C  \\
{C}^* & B \\
\end{array}\right).
$$
First we need to define $L_{k}$ on the words of length $2k-1$.  This will define the corner $C$ in the above matrix. 
Let us enumerate the words of length $k-2$ by $z_1, \ldots, z_t$ and the words of length $k-1$ by $w_1, \ldots, w_s$. 
By the definition of the tip-maximal functional we can find generators of $\ker A = \ker B_{L_{k-1}}$ of the form  $\sum_{i} \alpha_i^{(m)}w_i  - \sum_{j}\beta_j^{(m)} z_j$  ($m=1,\ldots,r$) such that elements 
$f_m= \sum_{i} \alpha_i^{(m)}w_i$ are linearly independent. 
Thus for every word $v\in S_{k-1}$ we have 
\begin{equation}\label{rel}
 \sum_{i} \alpha_i^{(m)} L_{k-1}(w_i v)  =  \sum_{j}\beta_j^{(m)} L_{k-1}(z_j v). 
\end{equation}

$S_{2k-1}$, clearly, coincides with  the set of   words $w_i w$ were  $\abs{w} = k$ and $i=1,\ldots, s$.   For each $w_i w$ we introduce a variable $x_{i,w}$ and for each $j=1,\ldots, t$ a constant $a_{j, w}= L_{k-1}(z_j w)$ ($z_j w \in S_{2k-2}$ hence $L_{k-1}(z_j w)$ is defined). And consider the system of linear equations 
\begin{equation} \label{myeq}
\sum_{i} \alpha_i^{(m)} x_{i,w}  =  \sum_{j}\beta_j^{(m)} a_{j, w}, \quad  (m=1, \ldots, r).
\end{equation}
Since the family  $w_1 w, \ldots, w_s w$ and the family $f_1, \ldots, f_r$ are both  linearly independent  the matrix $(\alpha_{i}^{(m)})_{i,m}$ has rank $r$. Thus  system~\eqref{myeq}  has a solution. Given a solution $x_{j,m}$ define 
$$
L_{k}(v) =\begin{cases} 
L_{k-1}(v) & \text{ for } v\in S_{k-1},\\
 x_{i,w}  & \text{ for }  v= w_i w \text{ with some } w \text{ such that }  \abs{w} = k.
 \end{cases}
$$
Thus $L_k$ is defined on $V_{2k-1}$. Hence the operator  $C \colon \Span( S_k\setminus S_{k-1}) \to V_{k-1}$ is defined by the matrix with $(v, w)$-entry equal to $L_k(v w^*)$ where $v\in V_{k-1}$ and $w\in S_k\setminus S_{k-1}$. 
By the definition of $x_{i,w}$ we immediately have that  
\begin{equation}\label{mainrel}
 \sum_{i} \alpha_i^{(m)} L_{k}(w_i w)  =  \sum_{j}\beta_j^{(m)} L_{k}(z_j w). 
\end{equation}  
for $w\in S_k\setminus S_{k-1}$. Since $L_{k}$ is extension of $L_{k-1}$ the condition~\eqref{mainrel} is satisfied also for all $w \in S_{k}$.   
Since for every generator $\sum_{i} \alpha_i^{(m)}w_i  - \sum_{j}\beta_j^{(m)} z_j$ of $\ker B_{L_{k-1}}$ equation \eqref{rel} implies \eqref{mainrel} we get that the subspace generated by  vectors 
$$(L_{k-1}(z_1 v),\ldots, L_{k-1}(z_s v), L_{k-1}(w_1 v), \ldots, L_{k-1}(w_r v)), $$ where $v\in V_{k-1}$, contains the  subspace generated by $$(L_{k}(z_1 u),\ldots, L_{k}(z_s u), L_{k}(w_1 u), \ldots, L_{k}(w_r u)), $$ where $u\in \Span(S_k\setminus S_{k-1}).$
 Thus $\Ran(C) \subseteq \Ran(A)$. 
 
 Since  the words $uv^*$ with $\abs{u} =\abs{v} = k$ are distinct for different pairs $(u,v)$ they  constitute a linearly independent family (which is $S_{2k}\setminus S_{2k-1}$). Thus we can define $L_k(uv^*) = B_{u,v}$ where $B= C^*(A|_{\Ran(A)})^{-1} C$.   By Lemma~\ref{lemF},  $L_k$ is flat and is positive semidefinite if $L_{k-1}$ is such.    
 
\end{proof}

As a corollary we have the following 
\begin{theorem} 
Every positive semidefinite truncated tip--maximal functional on a free $*$-algebra $\F$ has a positive semidefinite extension to  $\F$.  
\end{theorem}
\begin{proof} 
Any positive semidefinite truncated tip--maximal functional $L_{k}$  can be extended to a flat positive semidefinite $L_{k+1}$ by Lemma~\ref{lemRec} and then $L_{k+1}$ can be extended again to a positive semidefinite $L$ on $\F$ by Theorem~\ref{thmExt}.  

\end{proof}

\noindent {\bf Example.}
 
Consider the $*$-algebra $A =\Complex \seq{ x, x^* \mid x^2=0, x^{*2} = 0}$. Then 
$S_3 =\{x, x^*, x x^*, x^* x,$  $ x x^* x, x^* x x^*\}.$ 
By direct calculation we have that any hermitian $L_3$ with real coefficients  has a form 
$$
\left(\begin{array}{cccccc}
a_1 & 0 & 0 & a_3 & a_5 & 0 \\
0 & a_2 & a_4 & 0 & 0 & a_6 \\
0 & a_4 & a_5 & 0 & 0 & a_7 \\
a_3 & 0 & 0 & a_6 & a_8 & 0 \\
a_5 & 0 & 0 & a_8 & a_9 & 0 \\
0 & a_6 & a_7 & 0 & 0 & a_{10} \\
\end{array}\right)
$$ 

By Lemma~\ref{lemF}, $L_3$ is flat if and only if 

\begin{gather}
a_{9} = \frac{a_6^2 a_5 - 2 a_3 a_5 a_8 + a_1 a_8^2}{a_1 a_6 - a_3^2}, \label{eq1}\\ 
a_{10} = \frac{a_5 a_8^2 - 2 a_4 a_6 a_7 + a_2 a_7^2}{a_2 a_5 - a_4^2}. \label{eq2}
\end{gather} One can check that $L_2$ is positive definite iff $a_2>0$,  $a_6>0$, $a_2 a_5> a_4^2$ and $a_1 a_6 > a_3^2$. Thus any positive $L_2$ can be extended to a positive semidefinite flat $L_3$.

\noindent{\bf Example.}

Consider the free $*$-algebra $\F = \Complex \seq{x, x^*}$. Denote by $\B_+$ the set of words of the form $w w^*$ with $w$ is an arbitrary word in generators $x$  and $x^*$. It was shown in~\cite{Popovych:RepHoust} that there are faithful positive functionals on $\F$ (as well as on any Lance-Tapper $*$-algebra) of the form 
 $a\mapsto F(\Pi(a))$ where $\Pi \colon \F \to \F $ is linear extension of the map $\Pi \colon \B \to \B$ such that  
$$  \Pi(w) =\begin{cases}   
 w & \text{ if } w\in \B_+  \\
  0 & \text{ if } w \not\in \B_+
  \end{cases}
$$
and $F\colon \Complex \B_+ \to \Complex$ is a linear functional on subspace generated by $\B_+$. 

Assume that $L_4 \in V_8^*$ is such that 

\begin{equation}\label{poscond}
{L_{4}}|_{V_7} \text{ is of the form } F(\Pi|_{V_7}) \text{ for some functioal } F\in V_7^*
\end{equation}
 and has a block matrix decomposition $$\left(\begin{array}{cc}
A & C  \\
{C}^* & B \\
\end{array}\right)
$$ such that matrix $A$ is of the form  
$$ A=  \left( \begin{array}{cc}  
A_{11} & A_{12} \\
A_{12}^* & A_{22} 
\end{array}\right)
$$
with  
$A_{11} = I_{14}$, the identity $14\times 14$-matrix, and $A_{22}= \diag(2,3,1,1,2,3,0)$

$$
A_{12} = \left( \begin{array}{ccccccc} 
1 & 1 & 0 & 0 & 0 & 0 & 0 \\
0 & 0 & 0 & 0 & 1 & 1 & 0 \\
0 & 0 & 0 & 0 & 0 & 0 & 0\\
0 & 0 & 0 & 0 & 0 & 0 & 0\\
0 & 0 & 0 & 0 & 0 & 0 & 0\\
0 & 0 & 0 & 0 & 0 & 0 & 0\\
0 & 0 & 0 & 0 & 0 & 0 & 0\\
\end{array}\right)
$$
By easy but routine calculation one can check that the form on $V_3 \times V_3$ defined by the matrix $A$ comes from the functional  of the form $F(\Pi|_{V_6})$ and is positive semidefinite. Condition~\eqref{poscond} ensures that  the matrix $C$ is  completely determined by $A$. 

In the decomposition $V_3 = A(V_3) \oplus \ker A$ the subspace $A(V_3)$ is $$\Span \{ x, x^*, x^2, x x^*,  x^* x, x^{*2}, x^3, x^2 x^*, x x^* x, x x^{*2}, x^* x^2, x^* x x^*, x^{*2} x \}$$ and $\ker A  = \Span \{ x^{*3} \}.$
The matrix $A$ is chosen such that it is positive semidefinite tip--maximal and annulates $x^{3*}$.    The matrix $B$ and $L_4$ can be completely determined by Lemma~\ref{lemF}. 

The linear basis of the $\ker B_{L_4}$ is the following set $H$ 
\begin{gather*}
x^{*3}, x^4, x^3 x^* -x^2, x^2 x^* x- 2 x^2, x^2  x^{*2}, x x^* x^2, (x x^*)^2 -3 x x^*,  x x^{*2} x -x x^*, x x^{*3},\\  x^* x^3, x x^{*3}, x^*x^3, x^* x^2 x^* -x^* x, (x^* x)^2 - 2 x^* x, x^*x x^{*2}, x^{*2} x^2, x^{*2} x x^* - 3 x^{*2}, x^{*3} x, x^{*4} 
\end{gather*} 

A right Gr\"obner basis $\G$ of the right ideal $J$ generated by $\ker B_{L_4}$ is  
\begin{gather*}
x^{*3}, x^4, x^3 x^* -x^2, x^2 x^* x- 2 x^2, x^2  x^{*2}, x x^* x^2, (x x^*)^2 -3 x x^*,  x x^{*2} x -x x^*, x x^{*3},\\  x^* x^3,  x^* x^2 x^* -x^* x, (x^* x)^2 - 2 x^* x, x^*x x^{*2}, x^{*2} x^2, x^{*2} x x^* - 3 x^{*2}
\end{gather*} 
The only reductions we have made to compute  the above Gr\"obner basis starting from $H$ are the reductions of $x^{*3} x$ and  $ x^{*4}$ to zero. 

The positive definite form $\widetilde{B}_{L_4}$ is the form with a matrix $\widetilde{A}$ in the basis consisting of cosets of $\F / J$ with the representatives  $$x, x^*, x^2, x x^*, x^* x, x^{*2}, x^3, x^2 x^*, x x^* x, x x^{*2}, x^* x^2, x^* x x^*, x^{*2} x.$$ The matrix $\widetilde{A}$ is  obtained from the matrix $A_{11}$ by deleting the first row and first column.  Operators $R_x$ and $R_{x^*}$ of multiplications from the right by $x$ and $x^*$  define mutually adjoint operators $\widetilde{R}_x$ and $\widetilde{R}_{x^*}$ in the Hilbert space of right cosets $\F / J$ with inner product given by $\widetilde{A}$. Thus we get 13 dimensional representation $\pi$ of the free $*$-algebras $\F$. It can be checked that $\ker \pi \cap V_4$ is linearly generated by the following elelmets:
\begin{gather*}
x^3, -5 x^{*2}+2 x^* x x^{*2}+x^{*2} x x^*, x^* x^2 x^*-x  x^{*2} x.
\end{gather*}   
In particular,  $\pi$ defines a representation of Lance-Tapper $*$-algebra 
$$\Complex \seq{x, x^* \mid x^3=0, x^{*3} = 0}.$$

\end{document}